 \documentclass[final,1p,times]{elsarticle}

\usepackage[colorlinks=true]{hyperref}

\usepackage{amsthm,amsmath,amssymb,mathrsfs,amsfonts,graphicx,graphics,latexsym,exscale,cmmib57,dsfont,bbm,amscd,ulem}
\usepackage{color,soul}
\newtheorem{thm}{Theorem}[section]

\newtheorem{cor}[thm]{Corollary}
\newtheorem{lem}[thm]{Lemma}

\newcommand{\lam}[1]{\lambda_}
\newcommand{\bmat}[1]{\begin{bmatrix}#1\end{bmatrix}}
\newcommand{\bea}{\begin{eqnarray}}
\newcommand{\eea}{\end{eqnarray}}
\newcommand{\bean}{\begin{eqnarray*}}
\newcommand{\eean}{\end{eqnarray*}}
\newtheorem{formula}[thm]{Berger-Wang Formula}
\newtheorem{problem}[thm]{Problem}
\newtheorem{Barabanov}[thm]{Barabanov's norm theorem}
\newtheorem{defn*}{Definitions}

\newtheorem{Example}{Example}

\newcommand{\bcdot}{\boldsymbol{\cdot}}
\newcommand{\bS}{\boldsymbol{S}}

\newcommand{\set}[1]{\left\{#1\right\}}

\newcommand{\N}{{\mathbb N}}

\newcommand{\eps}{\varepsilon}

\journal{Linear Algebra and its Applications }
\begin{document}

\begin{frontmatter}
\title{The finite-step realizability of the joint spectral radius of a pair of $d\times d$ matrices one of which being rank-one\tnoteref{t1}}
\tnotetext[t1]{This work was supported in part by National Science Foundation of China (Grant Nos. 11071112 and 11071263) and in part by NSF 1021203 of the United States.}

\author[Dai]{Xiongping Dai}
\ead{xpdai@nju.edu.cn}
\author[Huang]{Yu Huang}
\ead{stshyu@mail.sysu.edu.cn}
\author[LX]{Jun Liu}
\ead{jliu@math.siu.edu}
\author[LX]{Mingqing Xiao}
\ead{mxiao@math.siu.edu}


\address[Dai]{Department of Mathematics, Nanjing University, Nanjing 210093, People's Republic of China}
\address[Huang]{Department of Mathematics, Zhongshan (Sun Yat-Sen) University, Guangzhou 510275, People's Republic of China}
\address[LX]{Department of Mathematics, Southern Illinois University, Carbondale, IL 62901-4408, USA}

\begin{abstract}
We study the finite-step realizability of the joint/generalized spectral radius of a pair of real $d\times d$ matrices $\{\mathrm{S}_1,\mathrm{S}_2\}$, one of which has rank $1$, where $2\le d<+\infty$.
Let $\rho(A)$ denote the spectral radius of a square matrix $A$. Then we prove that there always exists a finite-length word $({i}_1^*,\ldots,{i}_\ell^*)\in\{1,2\}^\ell$, for some finite $\ell\ge1$, such that
\begin{equation*}
\sqrt[\ell]{\rho(\mathrm{S}_{{i}_1^*}\cdots \mathrm{S}_{{i}_\ell^*})}=\sup_{n\ge1}\left\{\max_{(i_1,\ldots,i_n)\in\{1,2\}^n}\sqrt[n]{\rho(\mathrm{S}_{i_1}\cdots \mathrm{S}_{i_n})}\right\};
\end{equation*}
that is to say, there holds the spectral finiteness property for $\{\mathrm{S}_1,\mathrm{S}_2\}$. This implies that stability is algorithmically decidable for $\{\mathrm{S}_1,\mathrm{S}_2\}$.
\end{abstract}


\end{frontmatter}

\section{Introduction}\label{sec1}%
Let $\bS=\{\mathrm{S}_1,\ldots,\mathrm{S}_K\}\subset \mathbb{R}^{d\times d}$ be an arbitrary finite set of real $d$-by-$d$ matrices and
$\|\cdot\|$ a matrix norm on space $\mathbb{R}^{d\times d}$, where $2\le d<+\infty$ and $K\ge 2$. To capture the maximal growth rate of the trajectories generated by random products of matrices $\mathrm{S}_1,\ldots,\mathrm{S}_K$ in $\bS$, in 1960 \cite{RS} G.-C.~Rota and G.~Strang introduced the very important concept\,--\,{\it joint spectral radius of $\bS$}\,--\,by
\begin{equation*}
\hat{\rho}(\bS)=\lim_{n\to+\infty}\left\{\max_{(i_1,\ldots,i_n)\in\mathbb{K}^n}\sqrt[n]{\|\mathrm{S}_{i_1}\cdots \mathrm{S}_{i_n}\|}\right\}.
\end{equation*}
Here
$\mathbb{K}^n:=\overset{n\textrm{-time}}{\overbrace{\{1,\ldots,K\}\times\cdots\times\{1,\ldots,K\}}}$
stands for the set of all words $(i_1,\ldots,i_n)$ of finite-length $n$, composed by the letters $1,\ldots, K$, for any integer $n\ge1$. Let
\begin{equation*}
\varSigma_{\!K}^+=\{i_{\bcdot}\colon\mathbb{N}\rightarrow\mathbb{K}\},\quad \textrm{where }\mathbb{N}=\{1,2,\ldots\},
\end{equation*}
be the set of all one-sided infinite sequences (also called switching signals of $\bS$). Then we see, from N.~Barabanov~\cite{Bar} for example, that $\hat{\rho}(\bS)<1$ if and only if
\begin{equation*}
\|\mathrm{S}_{i_1}\cdots\mathrm{S}_{i_n}\|\to0\quad \textrm{as }n\to+\infty\qquad \forall i_{\bcdot}\in\varSigma_{\!K}^+.
\end{equation*}
In other words, $\hat{\rho}(\bS)<1$ if and only if the linear switched dynamical system, also write as $\bS$,
\begin{equation*}
x_n=x_0\cdot\mathrm{S}_{i_1}\cdots\mathrm{S}_{i_n},\quad x_0\in\mathbb{R}^d\textrm{ and }n\ge1,\qquad i_{\bcdot}\in\varSigma_{\!K}^+,
\end{equation*}
it is absolutely asymptotically stable, where the initial state $x_0\in\mathbb{R}^d$ is thought of as a row vector. In fact, from \cite{DHX-ERA} there follows
\begin{equation*}
\hat{\rho}(\bS)=\max_{i_{\bcdot}\in\varSigma_{\!K}^+}\left\{\limsup_{n\to+\infty}\sqrt[n]{\|\mathrm{S}_{i_1}\cdots \mathrm{S}_{i_n}\|}\right\}.
\end{equation*}
So, $\hat{\rho}(\bS)$ is a nonnegative real number which is independent of the norm $\|\cdot\|$ used here. It is a well-known fact that
the joint spectral radius $\hat{\rho}$ plays a critical role in a variety of applications such as switched dynamical systems \cite{Bar, Koz90, BW92, Gur95, Wirth02, Wirth05, Bar05, GWZ05, LD06-1, LD06-2, SWMW07, Jun09, DHX11, Dai-JDE}, differential equations \cite{Bar89, GZ, Dai-JDE}, coding theory \cite{MOS}, wavelets \cite{DL92-01, DL92, HS95, Mae98}, combinatorics~\cite{DST}, and so on.

Although $\hat{\rho}(\bS)$ is independent of the norm $\|\cdot\|$ used here, its approximation based on the above limit definition does rely upon an explicit choice of the norm $\|\cdot\|$ and has been a substantially interesting topic, for example, in \cite{Koz90, LW95, Gri96, Pro96, Pro97, TB97, Mae98, BN05, JPB08, PJ08, Koz09-1, Koz09-2, Koz10-DEDS, Koz10, Morris-am}. In general, computing $\hat{\rho}$
by definition cannot stop at some finite-time $n$, as shown by the single matrix
\begin{equation*}
A=\left[\begin{matrix}1&0\\ 1&1\end{matrix}\right]
\end{equation*}
where $\hat{\rho}(A)=\lim_{n\to+\infty}\sqrt[n]{\|A^n\|}=1$
by the classical Gel'fand spectral radius formula, however, there holds
$\sqrt[n]{\|A^n\|}>1$ for all $n\ge1$. For that reason in part, I.~Daubechies and J.~Lagarias in 1992~\cite{DL92-01}
defined the equally important concept -- {\it generalized spectral radius of $\bS$} -- by
\begin{equation*}
{\rho}(\bS)=\limsup_{n\to+\infty}\left\{\max_{(i_1,\ldots,i_n)\in\mathbb{K}^n}\sqrt[n]{\rho(\mathrm{S}_{i_1}\cdots \mathrm{S}_{i_n})}\right\},
\end{equation*}
where $\rho(A)$ stands for the usual spectral radius for any
matrix $A\in\mathbb{R}^{d\times d}$. And they conjectured there that a Gel'fand-type formula should hold for $\mathrm{S}$. This was proved by M.A.~Berger and Y.~Wang in 1992~\cite{BW92}, i.e., there holds the following Gel'fand-type formula.

\begin{formula}
${\rho}(\bS)=\hat{\rho}(\bS)$, for any bounded subset $\bS\subset\mathbb{R}^{d\times d}$.
\end{formula}

Because of its importance, this Gel'fand-type spectral-radius formula has been reproved by using different interesting approaches,
for example, in~\cite{El, SWP, CZ00, Bochi, Dai-JMAA}. According to this formula, the computation of ${\rho}(\bS)$ becomes an important subject at once, which leads to the following significant problem motivated by
\begin{equation*}
\rho(A)=\sqrt[n]{\rho(A^n)}
\end{equation*}
for any square matrix $A$.

\begin{problem}[Spectral Finiteness Property]\label{problem1.2}
Does there exist any word $(i_1,\ldots,i_n)$ of finite-length $n\ge1$ such that
\begin{equation*}
{\rho}(\bS)=\sqrt[n]{{\rho}(\mathrm{S}_{i_1}\cdots \mathrm{S}_{i_n})},
\end{equation*}
for any $\bS=\{\mathrm{S}_1,\ldots,\mathrm{S}_K\}\subset\mathbb{R}^{d\times d}$?
\end{problem}

This spectral finiteness property means that ${\rho}(\bS)$ is computationally efficient. It was conjectured, respectively, by E.S.~Pyatnitski\v{i}~\cite{PR91} for its continuous-time version, I.~Daubechies and J.~Lagarias in \cite{DL92-01}, L.~Gurvits in~\cite{Gur95}, and by J.~Lagarias and Y.~Wang in \cite{LW95}. If this were true for $\mathrm{S}$, then from the Berger-Wang formula, it follows that we would realize efficiently the joint/generalized spectral radius $\hat{\rho}(\bS)$ only by computation of finite steps, and the interests would arise from its connection with the stability question for $\bS$.

Unfortunately, this important ``spectral finiteness conjecture" has been disproved
by T.~Bousch and J.~Mairesse in \cite{BM} using measure-theoretical ideas, also respectively by V.~Blondel \textit{et al.}
in \cite{BTV} exploiting combinatorial properties of permutations of products of positive matrices, and by V.~Kozyakin~\cite{Koz05, Koz07} employing the theory of dynamical systems, all offered the existence of counterexamples in the case where $d=2$ and $K=2$. Moreover, an explicit expression for such a counterexample has been found in the recent work of K.~Hare \textit{et al.}~\cite{HMST}.

Although the finiteness conjecture fails to exist, the idea of Problem~\ref{problem1.2} is
still to be very attractive and important due to developing efficient algorithms because the
computation of the joint spectral radius $\hat{\rho}$ must be implemented in
finite arithmetic. Some conjectures in special case still keep open, for example, in M.~Maesumi~\cite{Mae96} and R.~Jungers and V.~Blondel~\cite{JB08}. Many positive efforts have been made and studies show that spectral finiteness property may be true in a number of interesting cases, for example, see \cite{Gur95, LW95, GWZ05, Theys, BJP06, JB08, GZ01, GZ08, GZ09, JP09, CGSZ10, Morris10}, including the case were the matrices $\mathrm{S}_1,\ldots,\mathrm{S}_K$ are symmetric, or if the Lie algebra associated with the set of matrices is solvable~\cite[Corollary 6.19]{Theys}; in this case $\rho(\bS)=\max_{1\le i\le K}\rho(\mathrm{S}_i)$, see \cite{Gur95, LHM, JB08}. Particularly, in A.~Cicone \textit{et al.} \cite{CGSZ10} it was proved, based on R.~Jungers and V.~Blondel~\cite{JB08} which is for all pairs of $2\times 2$ binary matrices, that
every pairs of $2\times 2$ sign-matrices $\mathrm{S}_1,\mathrm{S}_2$ have the spectral finiteness property described in Problem~\ref{problem1.2}.

In the present paper, based on the important work of Barabanov~\cite{Bar}, we will prove, mathematically in Section~\ref{sec2} and numerically in Section~\ref{sec3}, the following finiteness result.

\begin{thm}\label{thm1.3}
Let $2\le d<+\infty$ and $\bS=\{\mathrm{S}_1, \mathrm{S}_2\}$ be an arbitrary pair of real $d\times d$ matrices.
If one of $\mathrm{S}_1$ and $\mathrm{S}_2$ has rank $1$, then $\bS$ has the spectral finiteness property.
\end{thm}

This means, from \cite[Proposition~1]{JB08}, that stability is algorithmically decidable, for every pairs of real $d\times d$ matrices $\mathrm{S}_1,\mathrm{S}_2$ if one of which has rank $1$. If, in addition, $\bS$ is irreducible, then $\bS$ possesses the rank one property introduced by I.D.~Morris in~\cite{Morris10}. However, the counterexample of Hare \textit{et al.} \cite{HMST} shows that Morris's rank one property is neither necessary nor sufficient for the finiteness property. So, our rank $1$ condition described in Theorem~\ref{thm1.3} is substantial for our statement.

By $\bS^+$, it means for the multiplicative semigroup generated by $\mathrm{S}_1, \mathrm{S}_2$, i.e,
\begin{equation*}
\bS^+={\bigsqcup}_{n\geq1}\left\{\mathrm{S}_{i_1}\cdots\mathrm{S}_{i_n}\,|\,(i_1,\ldots,i_n)\in\mathbb{K}^n\right\}\quad \textrm{where }\mathbb{K}=\{1,2\}.
\end{equation*}
Here $\bigsqcup$ means the disjoint union.
To mathematically prove Theorem~\ref{thm1.3}, our approach is to consider its equivalent statement formulated as follows:

\begin{thm}\label{thm1.4}
Let $\bS=\{\mathrm{S}_1, \mathrm{S}_2\}$ be an arbitrary pair of real $d\times d$ matrices, one of which has rank $1$.
If $\rho(A)<1$ for all $A\in\bS^+$, then $\rho(\bS)<1$; namely, the induced switched dynamics $\bS$ is absolutely asymptotically and exponentially stable.
\end{thm}

This stability was first conjectured by E.S.~Pyatnitski\v{i} in 1980s, see, e.g.,~\cite{PR91, Gur95, SWMW07} and it has been the subject of substantial recent research interest, for example, in \cite{Gur95, SWP, SWMW07, DHX11, Dai-JDE}.

This paper is organized as follows. In Section~\ref{sec2}, we will provide the proof of our main results. Moreover, we also give an explicit formula for the computation of the generalized spectral radius. Several examples are provided in Section~\ref{sec3} to illustrate the theoretical outcomes. The paper ends with concluding remarks in Section~\ref{sec4}.
\section{Finiteness property of a pair of $d\times d$ matrices}\label{sec2}
This section is devoted to mathematically proving our main results, Theorems~\ref{thm1.3} and \ref{thm1.4}.

To prove Theorem~\ref{thm1.3}, we first prove the following absolute stability theorem, which is important not just to our spectral finiteness theorem, but also to the stabilizability of switched dynamical systems \cite{DHX11, SWMW07}.

\begin{thm}\label{thm2.1}
Let $\bS=\{\mathrm{S}_1,\mathrm{S}_2\}\subset\mathbb{R}^{d\times d}$, $2\le d<+\infty$, be periodically switched stable; that is to say,
\begin{equation*}
\rho(A)<1\quad\forall A\in\bS^+.
\end{equation*}
Then, if one of $\mathrm{S}_1, \mathrm{S}_2$ is of rank $1$, $\bS$ is absolutely exponentially stable, i.e.,
\begin{equation*}
\|\mathrm{S}_{i_1}\cdots \mathrm{S}_{i_n}\|\xrightarrow[]{\textrm{exponentially fast}}0\quad \textrm{as }n\to+\infty,\qquad\textrm{i.e., }\limsup_{n\to+\infty}\frac{1}{n}\log\|\mathrm{S}_{i_1}\cdots \mathrm{S}_{i_n}\|<0,
\end{equation*}
for all switching signals $i_{\bcdot}\colon \N\rightarrow\{1,2\}$.
\end{thm}

Here $\bS^+$ is the multiplicative semigroup generated by $\bS$ as described in Section~\ref{sec1}. Recall that $\bS$ is said to be irreducible, provided that there is no common, nontrivial and proper invariant linear subspaces of $\mathbb{R}^d$, for $\mathrm{S}_1, \mathrm{S}_2$.

The following result holds trivially by induction on $d$ together with the Berger-Wang formula, which is a standard result in the theory of linear algebras.

\begin{lem}[{See, e.g.,\cite{Bar, BW92, Dai-JMAA}}]\label{lem2.2}
For any $\bS=\{\mathrm{S}_1,\mathrm{S}_2\}\subset\mathbb{R}^{d\times d}$, there exists a nonsingular matrix $\mathrm{P}\in\mathbb{R}^{d\times d}$ and $r$ positive integers $d_1,\ldots, d_r$ with $d_1+\cdots+d_r=d$ such that
\begin{equation*}
\mathrm{P}\mathrm{S}_i\mathrm{P}^{-1}=\left[\begin{matrix}\mathrm{S}_i^{(1,1)}&\mathbf{0}_{d_1\times d_2}&\cdots&\mathbf{0}_{d_1\times d_r}\\ \mathrm{S}_i^{(2,1)}&\mathrm{S}_i^{(2,2)}&\cdots&\mathbf{0}_{d_2\times d_r}\\ \vdots&\vdots&\ddots&\vdots\\ \mathrm{S}_i^{(r,1)}&\mathrm{S}_i^{(r,2)}&\cdots&\mathrm{S}_i^{(r,r)}\end{matrix}\right]\quad (i=1,2),
\end{equation*}
where $\bS^{(k)}:=\left\{\mathrm{S}_1^{(k,k)}, \mathrm{S}_2^{(k,k)}\right\}\subset\mathbb{R}^{d_k\times d_k}$ is irreducible for each $1\le k\le r$, such that
\begin{equation*}
\max\left\{\rho(\bS^{(k)})\colon 1\le k\le r\right\}=\rho(\bS)=\hat{\rho}(\bS)=\max\left\{\hat{\rho}(\bS^{(k)})\colon 1\le k\le r\right\}.
\end{equation*}
\end{lem}

When $\bS$ is itself irreducible, $r=1$ in Lemma~\ref{lem2.2}.

The following important theorem, due to N.~Barabanov, is extremely valuable to the proof of Theorem~\ref{thm2.1}.

\begin{Barabanov}[{See~\cite{Bar}, also \cite{Wirth02, Dai-JMAA}}]\label{thm2.3}
If $\bS=\{\mathrm{S}_1, \mathrm{S}_2\}\subset\mathbb{R}^{d\times d}$ is irreducible, then there is a vector norm $|\pmb{|}\cdot\pmb{|}|_*$ on $\mathbb{R}^d$ such that there hold the following two statements.
\begin{enumerate}
\item[$\mathrm{(1)}$] $\hat{\rho}(\bS)=\max\left\{\sqrt[n]{|\pmb{|}\mathrm{S}_{i_1}\cdots\mathrm{S}_{i_n}\pmb{|}|_*}\colon (i_1,\ldots,i_n)\in\mathbb{K}^n\right\}$ for all $n\ge1$.

\item[$\mathrm{(2)}$] To any $\hat{x}\in\mathbb{R}^d$, there corresponds an infinite sequence, say $i_{\bcdot}(\hat{x})\colon\mathbb{N}\rightarrow\mathbb{K}$,
satisfying that
$|\pmb{|}\hat{x}\cdot \mathrm{S}_{i_1(\hat{x})}\cdots \mathrm{S}_{i_n(\hat{x})}\pmb{|}|_*=\hat{\rho}(\bS)^n|\pmb{|}\hat{x}\pmb{|}|_*$
for all $n\ge1$.
\end{enumerate}
Here $\mathbb{K}=\{1,2\}$ and the matrix norm $|\pmb{|}\cdot\pmb{|}|_*$ on $\mathbb{R}^{d\times d}$ is naturally induced by the norm $|\pmb{|}\cdot\pmb{|}|_*$ on $\mathbb{R}^d$.
\end{Barabanov}

Using the Berger-Wang formula, Lemma~\ref{lem2.2} and Barabanov's norm theorem, we now can prove Theorem~\ref{thm2.1}.

\begin{proof}[\textbf{Proof of Theorem~\ref{thm2.1}}]
Let $\bS=\{\mathrm{S}_1, \mathrm{S}_2\}\subset\mathbb{R}^{d\times d}$ with $\mathrm{rank}(\mathrm{S}_2)=1$. Then according to Lemma~\ref{lem2.2}, there is no loss of generality is assuming that $\bS$ is irreducible. Since $\bS$ is periodically switched stable, we have $\hat{\rho}(\bS)\le1$ by the definition of $\rho(\bS)$ and the Berger-Wang formula. Therefore, from Barabanov's theorem, it follows that there exists a vector norm $|\pmb{|}\cdot\pmb{|}|_*$ on $\mathbb{R}^d$, which induces a matrix norm, write also $|\pmb{|}\cdot\pmb{|}|_*$, on $\mathbb{R}^{d\times d}$ such that
\begin{equation*}
|\pmb{|}\mathrm{S}_1\pmb{|}|_*\le1\quad\textrm{and}\quad|\pmb{|}\mathrm{S}_2\pmb{|}|_*\le1.
\end{equation*}
We simply write
\begin{equation*}
\mathrm{S}_1=\left[a_{ij}\right]_{d\times d}\quad \textrm{and}\quad \mathrm{S}_1^{\ell}=\overset{\ell\textrm{-time}}{\overbrace{\mathrm{S}_1\cdots\mathrm{S}_1}}=\left[a_{ij}^{(\ell)}\right]_{d\times d}\; \forall \ell\ge1.
\end{equation*}
As $\mathrm{S}_2$ is of rank $1$, it follows, from the Jordan canonical form theorem, that there is no loss of generality in assuming that
\begin{equation*}
\mathrm{S}_2=\mathrm{B}_1:=\left[\begin{array}{ll}\lambda&\mathbf{0}_{1\times (d-1)}\\ \mathbf{0}_{(d-1)\times 1}&\mathbf{0}_{(d-1)\times(d-1)}\end{array}\right] \quad \textrm{where } 0<|\lambda|<1
\end{equation*}
or
\begin{equation*}
\mathrm{S}_2=\mathrm{B}_2:=\left[\begin{matrix}0&0&0&\cdots&0\\ 1&0&0&\cdots&0\\0&0&0&\cdots&0\\\vdots&\vdots&\vdots&\ddots&\vdots\\0&0&0&\cdots&0\end{matrix}\right].
\end{equation*}
As $\{\mathrm{S}_1,\mathrm{S}_2\}$ is periodically switched stable, it follows from the classical Gel'fand spectral radius formula that
\begin{equation*}
|\pmb{|}\mathrm{S}_i^n\pmb{|}|_*\xrightarrow[]{\textrm{exponentially fast}}0\quad \textrm{as }n\to+\infty, \qquad\textrm{i.e., }\log\rho(\mathrm{S}_i)=\lim_{n\to+\infty}\frac{1}{n}\log|\pmb{|}\mathrm{S}_i^n\pmb{|}|_*<0,
\end{equation*}
for both $i=1$ and $2$.

Let $\mathbb{K}=\{1,2\}$. Next, we will prove the statement of Theorem~\ref{thm2.1} in the cases $\mathrm{S}_2=\mathrm{B}_1$ and $\mathrm{S}_2=\mathrm{B}_2$, respectively.

\textbf{Case I:} Let $\mathrm{S}_2=\mathrm{B}_1$. Note that in this case, for any finite-length word of the form
\begin{equation*}
w=(\mathfrak{i}_1,\ldots,\mathfrak{i}_\ell,\mathfrak{i}_{\ell+1},\ldots,\mathfrak{i}_{\ell+m})=(\stackrel{\ell\textrm{-time}}{\overbrace{1,\ldots,1}}, \stackrel{m\textrm{-time}}{\overbrace{2,\ldots,2}})\in\mathbb{K}^{\ell+m},
\end{equation*}
there holds
\begin{equation*}
\mathrm{S}(w):=\mathrm{S}_{\mathfrak{i}_1}\cdots \mathrm{S}_{\mathfrak{i}_\ell}\mathrm{S}_{\mathfrak{i}_{\ell+1}}\cdots\mathrm{S}_{\mathfrak{i}_{\ell+m}}=\mathrm{S}_1^\ell \mathrm{B}_1^m=\lambda^m\left[\begin{matrix}a_{11}^{(\ell)}&0&\cdots&0\\ a_{21}^{(\ell)}&0&\cdots&0\\\vdots&\vdots&\ddots&\vdots\\a_{d1}^{(\ell)}&0&\cdots&0\end{matrix}\right],
\end{equation*}
for any $\ell\ge1$ and $m\ge1$. Since $\{\mathrm{S}_1,\mathrm{S}_2\}$ is periodically switched stable, $\mathrm{S}(w)$ is exponentially stable and so it holds from the classical Gel'fand formula that
\begin{equation*}
\rho(\mathrm{S}(w))=|\lambda^m a_{11}^{(\ell)}|<1
\end{equation*}
for all words $w=(\stackrel{\ell\textrm{-time}}{\overbrace{1,\ldots,1}}, \stackrel{m\textrm{-time}}{\overbrace{2,\ldots,2}})$, for all $\ell\ge1$ and $m\ge1$.

Let $i_{\bcdot}\colon\N\rightarrow\mathbb{K}$ be an arbitrary switching signal. If to any $N\ge1$ there is some $n\ge N$ so that the infinite-length sequence $i_{\bcdot}=(i_1,i_2,\ldots)$ contains at least one of the following two sub-words of finite-length $n$
\begin{equation*}
(\stackrel{n\textrm{-time}}{\overbrace{1,\ldots,1}})\quad \textrm{and}\quad (\stackrel{n\textrm{-time}}{\overbrace{2,\ldots,2}}),
\end{equation*}
then it holds that
\begin{equation*}
|\pmb{|}\mathrm{S}_{i_1}\cdots \mathrm{S}_{i_n}\pmb{|}|_*\to0\quad \textrm{as }n\to+\infty.
\end{equation*}
Hence, we only need to consider the following special case:
\begin{equation*}
i_{\bcdot}=(\underset{w_1}{\uwave{\overset{\ell_1\textrm{-time}}{\overbrace{1,\ldots,1}}, \overset{m_1\textrm{-time}}{\overbrace{2,\ldots,2}}}},\underset{w_2}{\uwave{\overset{\ell_2\textrm{-time}}{\overbrace{1,\ldots,1}}, \overset{m_2\textrm{-time}}{\overbrace{2,\ldots,2}}}},\underset{\ldots}{\underset{{}}{\ldots}}, \underset{w_n}{\uwave{\overset{\ell_n\textrm{-time}}{\overbrace{1,\ldots,1}}, \overset{m_n\textrm{-time}}{\overbrace{2,\ldots,2}}}}, \ldots)
\end{equation*}
where $1\le \ell_n\le L$ and $1\le m_n\le M$ for all $n\ge1$, for some two positive integers $L\ge1$ and $M\ge1$. Therefore, there exists a positive constant \begin{equation*}
\gamma=\gamma(L,M)<1
\end{equation*}
such that
\begin{equation*}
|\lambda^{m_n} a_{11}^{(\ell_n)}|=\rho(\mathrm{S}(w_n))\le\gamma\quad\forall n\ge1.
\end{equation*}
Notice here that for the given special switching signal $i_{\bcdot}\colon\mathbb{N}\rightarrow\mathbb{K}$, $\gamma$ is independent of the extremal norm $|\pmb{|}\cdot\pmb{|}|_*$ of $\mathrm{S}$ used here.
From the fact that
\begin{equation*}
\begin{split}
\limsup_{n\to+\infty}\frac{1}{n}\log|\pmb{|}\mathrm{S}_{i_1}\cdots \mathrm{S}_{i_n}\pmb{|}|_*&=\limsup_{n\to+\infty}\frac{1}{J_n}\log|\pmb{|}\mathrm{S}_{i_1}\cdots \mathrm{S}_{i_{J_n}}\pmb{|}|_*\quad \textrm{where }J_n=\sum_{k=1}^n(\ell_k+m_k)\\
&=\limsup_{n\to+\infty}\frac{1}{\sum_{k=1}^n(\ell_k+m_k)}\log{\prod}_{k=1}^n|\lambda^{m_k} a_{11}^{(\ell_k)}|\\
&\le\frac{1}{L+M}\log\gamma<0
\end{split}
\end{equation*}
by \cite[Theorem~2.1]{Dai-sicon} and the triangularity of $\mathrm{S}(w_n)$,
it follows at once that
\begin{equation*}
|\pmb{|}\mathrm{S}_{i_1}\cdots \mathrm{S}_{i_n}\pmb{|}|_*\to0\quad \textrm{as }n\to+\infty.
\end{equation*}

Since the switching signal $i_{\bcdot}\colon\mathbb{N}\rightarrow\mathbb{K}$ is arbitrary here, this proves that $\{\mathrm{S}_1,\mathrm{S}_2\}$ is absolutely asymptotically stable.

\textbf{Case (II):} Let $\mathrm{S}_2=\mathrm{B}_2$. Noting that
\begin{equation*}
\mathrm{S}_1^\ell \mathrm{S}_2=\left[\begin{matrix}a_{12}^{(\ell)}&0&\cdots&0\\ a_{22}^{(\ell)}&0&\cdots&0\\\vdots&\vdots&\ddots&\vdots\\a_{d2}^{(\ell)}&0&\cdots&0\end{matrix}\right]\quad\forall \ell\ge1\quad\textrm{and}\quad \mathrm{S}_2^m=\mathbf{0}_{d\times d}\quad\forall m\ge2,
\end{equation*}
we can prove, by an argument similar to that of the case \textbf{(I)}, that $\{\mathrm{S}_1,\mathrm{S}_2\}$ is also absolutely asymptotically stable in this case.

Now combining the cases \textbf{(I)} and \textbf{(II)}, we see that
\begin{equation*}
|\pmb{|}\mathrm{S}_{i_1}\cdots\mathrm{S}_{i_n}\pmb{|}|_*\to0\quad \textrm{as }n\to+\infty
\end{equation*}
for all switching signals $i_{\bcdot}\colon\mathbb{N}\rightarrow\mathbb{K}$. Then, the statement of Theorem~\ref{thm2.1} follows immediately from the Fenichel uniformity theorem proven in~\cite{Fen}.

This completes the proof of Theorem~\ref{thm2.1}.
\end{proof}

If there is no the assumption of rank $1$ in the above Theorem~\ref{thm2.1}, then we can only guarantee that $\bS$ is exponentially stable almost surely in terms of some special probabilities from \cite{DHX11, Dai-JDE}.

As a result of Theorem~\ref{thm2.1}, we can obtain the following finiteness property.

\begin{thm}\label{thm2.4}
Let $\bS=\{\mathrm{S}_1,\mathrm{S}_2\}\subset\mathbb{R}^{d\times d}$, where $2\le d<+\infty$. If one of $\mathrm{S}_1,\mathrm{S}_2$ is of rank $1$, then $\bS$ has the spectral finiteness property; that is, one can found some finite $n\ge1$ such that
\begin{equation*}
\rho(\bS)=\max\left\{\sqrt[n]{{\rho(\mathrm{S}_{i_1}\cdots\mathrm{S}_{i_n})}}\colon(i_1,\ldots,i_n)\in\mathbb{K}^n\right\}.
\end{equation*}
Here $\mathbb{K}=\{1,2\}$.
\end{thm}

\begin{proof}
There is no loss of generality in assuming $\rho(\bS)=1$, by normalization of $\bS$ if necessary. Suppose, by contradiction, that
\begin{equation*}
\rho(A)<1\quad\forall A\in\bS^+.
\end{equation*}
Then from Theorem~\ref{thm2.1}, it follows that the switched dynamics induced by $\bS$ is absolutely exponentially stable. Thus $\hat{\rho}(\bS)<1$ from \cite{Bar} for example, and further $\rho(\bS)<1$ from the Berger-Wang formula \cite{BW92}. It is a contradiction to the assumption of $\rho(\bS)=1$.

This thus ends the proof of Theorem~\ref{thm2.4}.
\end{proof}

As a consequence of Theorem~\ref{thm2.4}, we can conclude the following result, which means that stability is algorithmically decidable for every pairs of real $d\times d$ matrices $\mathrm{S}_1,\mathrm{S}_2$ one of which has rank $1$.

\begin{cor}\label{cor2.5}
Denote $\mathbb{Z}_+=\{0, 1, 2, \dots\}$. For every pairs of real $d\times d$ matrices $\mathrm{S}_1,\mathrm{S}_2$ with $\mathrm{rank}(\mathrm{S}_2)=1$, we have
\begin{equation*}
\rho(\bS)=\max_{\ell, m\in \mathbb{Z}_+}\sqrt[\ell+m]{\rho(\mathrm{S}_1^{\ell}\mathrm{S}_2^m)}.
\end{equation*}
More specifically, we have
\begin{itemize}
\item if $\rho(\mathrm{S}_2)= 0$, then
\begin{equation*}
\rho(\bS)=\max\left\{\max_{\ell\in \N}\sqrt[\ell+1]{\rho(\mathrm{S}_1^{\ell}\mathrm{S}_2)}, \quad \rho(\mathrm{S}_1)\right\}
\end{equation*}
\item if $\rho(\mathrm{S}_2)\ne 0$, then
\begin{equation*}
\rho(\bS)=\max_{\ell, m\in \mathbb{Z}_+}\sqrt[\ell+m]{\rho(\mathrm{S}_1^{\ell}\mathrm{S}_2^m)}
\end{equation*}
\end{itemize}
\end{cor}

\begin{proof}
Without loss of generality, we may assume
\begin{equation*}
\mathrm{S}_2=\mathrm{B}_1:=\left[\begin{array}{ll}\lambda&\mathbf{0}_{1\times (d-1)}\\ \mathbf{0}_{(d-1)\times 1}&\mathbf{0}_{(d-1)\times(d-1)}\end{array}\right]
\end{equation*}
or
\begin{equation*}
\mathrm{S}_2=\mathrm{B}_2:=\left[\begin{matrix}0&0&0&\cdots&0\\ 1&0&0&\cdots&0\\0&0&0&\cdots&0\\\vdots&\vdots&\vdots&\ddots&\vdots\\0&0&0&\cdots&0\end{matrix}\right].
\end{equation*}
Similar to the previous proof of Theorem~\ref{thm2.1}, the possible optimal sequences should have the form
\begin{equation*}
i_{\bcdot}=(\underset{w_1}{\uwave{\overset{\ell_1\textrm{-time}}{\overbrace{1,\ldots,1}}, \overset{m_1\textrm{-time}}{\overbrace{2,\ldots,2}}}},\underset{w_2}{\uwave{\overset{\ell_2\textrm{-time}}{\overbrace{1,\ldots,1}}, \overset{m_2\textrm{-time}}{\overbrace{2,\ldots,2}}}},\underset{\ldots}{\underset{{}}{\ldots}}, \underset{w_n}{\uwave{\overset{\ell_n\textrm{-time}}{\overbrace{1,\ldots,1}}, \overset{m_n\textrm{-time}}{\overbrace{2,\ldots,2}}}}, \ldots),
\end{equation*}
by noting that
\begin{align*}
&\rho(\mathrm{B}_1^\ell\mathrm{S}(w_1w_2\cdots w_n))=\rho(\mathrm{S}(w_1^\prime w_2\cdots w_n))\quad \textrm{where }w_1^\prime=(\uwave{\overset{\ell_1\textrm{-time}}{\overbrace{1,\ldots,1}}, \overset{(m_1+\ell)\textrm{-time}}{\overbrace{2,\ldots,2}}})\\
\intertext{and}&\rho(\mathrm{B}_2^\ell\mathrm{S}(w_1w_2\cdots w_n))=0
\end{align*}
for all $\ell\ge1$ and $n\ge1$.
Denote
\begin{equation*}
i_{\bcdot}(n)=w_1w_2\cdots w_n
\end{equation*}
for any $n\ge1$.

If $\mathrm{S}_2=\mathrm{B}_1$ we then have
\begin{equation*}
\rho(\mathrm{S}(i_{\bcdot}(n)))=\prod_{k=1}^n\rho(\mathrm{S}_2)^{m_k}a_{11}^{(\ell_k)}
\end{equation*}
which yields a maximum when $w_1=w_2=\cdots =w_n$. In this case
\begin{equation*}\rho(\mathrm{S}(i_{\bcdot}(n))) =\rho({\mathrm{S}(w_1)}^n)=\rho(S(w_1))^n=\rho(\mathrm{S}_1^{\ell}\mathrm{S}_2^m)^n.\end{equation*}
Now if we let
\begin{equation*}
\alpha=\sup_{\ell, m\in \N}\sqrt[\ell+m]{\rho(\mathrm{S}_1^{\ell}\mathrm{S}_2^m)},
\end{equation*}
then we have
\begin{equation*}
\sqrt[|i_{\bcdot}(n)|]{\rho(\mathrm{S}(i_{\bcdot}(n)))}\leq \alpha.
\end{equation*}
This gives  $\rho(\bS)\le \alpha$. On the other hand, we know that
\begin{equation*}
\alpha=\sup_{\ell, m\in \N}\sqrt[\ell+m]{\rho(\mathrm{S}_1^{\ell}\mathrm{S}_2^m)}\le \rho(\bS)
\end{equation*}
This leads to $\sup_{\ell, m\in \N}\sqrt[\ell+m]{\rho(\mathrm{S}_1^{\ell}\mathrm{S}_2^m)}=\rho(\bS)$ and so $\max_{\ell, m\in \N}\sqrt[\ell+m]{\rho(\mathrm{S}_1^{\ell}\mathrm{S}_2^m)}=\rho(\bS)$ from Theorem~\ref{thm2.4}.

If $\mathrm{S}_2=\mathrm{B}_2$, the possible optimal sequence is given by
\begin{equation*}
i_{\bcdot}=(\underset{w_1}{\uwave{\overset{\ell_1\textrm{-time}}{\overbrace{1,\ldots,1}}, \overset{m_1\textrm{-time}}{\overbrace{2,\ldots,2}}}},\underset{w_2}{\uwave{\overset{\ell_2\textrm{-time}}{\overbrace{1,\ldots,1}}, \overset{m_2\textrm{-time}}{\overbrace{2,\ldots,2}}}},\underset{\ldots}{\underset{{}}{\ldots}}, \underset{w_n}{\uwave{\overset{\ell_n\textrm{-time}}{\overbrace{1,\ldots,1}}, \overset{m_n\textrm{-time}}{\overbrace{2,\ldots,2}}}}, \ldots)
\end{equation*}
with $m_i\equiv 1$. This corresponds the previous case by letting $m_1=m_2=\cdots=1$.

Thus, the proof of Corollary~\ref{cor2.5} is completed.
\end{proof}
\section{Illustrated Examples}\label{sec3}%
In this section we provide several examples to illustrate our theoretical outcomes proved in Section~\ref{sec2}. We here point out that it is unnecessary to transform the rank-one matrix $\mathrm{S}_2$
to its Jordan canonical form during practical calculation,
since the corresponding optimal sequence is invariant under similarity transformation.
Now, let us carry on the above analysis on the following examples.

\begin{Example}[See \cite{JB08}]\label{example1}
Let $\bS=\{\mathrm{S}_1,\mathrm{S}_2\}$, where
\begin{equation*}
\mathrm{S}_1=\bmat{1& 0\\1& 1},\quad\mathrm{S}_2=\bmat{0& 1\\0& 0}.
\end{equation*}
\end{Example}

\noindent Since
\begin{equation*}
\mathrm{S}_1^{\ell}=\bmat{1&0\\ \ell& 1},
\end{equation*}
 we have
\begin{equation*}
\mathrm{S}_1^{\ell}\mathrm{S}_2 =\bmat{1& 1 \\ 0& \ell}.
\end{equation*}
Then
\begin{equation*}
\rho(\mathrm{S}_1^{\ell}\mathrm{S}_2)=\ell.\end{equation*}
Hence
\begin{equation*}
\rho(\bS)=\max_{\ell\in \N} \sqrt[\ell+1]{\ell}=\sqrt[5]{4}.
\end{equation*}
This yields $\rho(\bS)=\sqrt[5]{4}$ and the corresponding optimal sequence is $\mathrm{S}_1^4\mathrm{S}_2$.

\begin{Example}\label{example2}
\begin{equation*}
\bS=\set{\mathrm{S}_1=\bmat{1& \frac{1}{\sqrt{2}}\\0& 1},\mathrm{S}_2=\bmat{1& \frac{\sqrt{3}}{2}\\-1& -\frac{\sqrt{3}}{2}}}.
\end{equation*}
\end{Example}

\noindent Notice that
\begin{equation*}
\mathrm{S}_1^{\ell}=\bmat{1&\frac{\ell}{\sqrt{2}}\\0 & 1},\
\qquad \mathrm{S}_2^m=\left(1-\frac{\sqrt{3}}{2}\right)^{m-1}\bmat{1& \frac{\sqrt{3}}{2}\\-1& -\frac{\sqrt{3}}{2}}.
\end{equation*}
and
\begin{equation*}
\mathrm{S}_1^{\ell}\mathrm{S}_2^m=\left(1-\frac{\sqrt{3}}{2}\right)^{m-1}\bmat{1&\frac{\ell}{\sqrt{2}}\\0 & 1}
\bmat{1\\-1}\bmat{1& \frac{\sqrt{3}}{2}}.
\end{equation*}
Thus we have
\begin{equation*}
\rho(\mathrm{S}_1^{\ell}\mathrm{S}_2^m)=\left(\frac{\ell}{\sqrt{2}}+\frac{\sqrt{3}}{2}-1\right)\left(1-\frac{\sqrt{3}}{2}\right)^{m-1}.
\end{equation*}
Hence
\bean
\rho(\bS)&=&\max_{\ell, m\in \N}
\sqrt[\ell+m]{\left(\frac{\ell}{\sqrt{2}}+\frac{\sqrt{3}}{2}-1\right)\left(1-\frac{\sqrt{3}}{2}\right)^{m-1}}\\
&=&\sqrt[6]{\frac{5}{\sqrt{2}}+\frac{\sqrt{3}}{2}-1}\\
&\approx& 1.226346> \max\left\{\rho(\mathrm{S}_1),\rho(\mathrm{S}_2)\right\},
\eean
where the maximum is attained at $(\ell,m)=(5,1)$ with the optimal sequence $\mathrm{S}_1^5 \mathrm{S}_2$.

\begin{Example}\label{example3}
\begin{equation*}
\bS=\set{\mathrm{S}_1=\bmat{1& \frac{1}{\sqrt{2}}\\0& 1},\mathrm{S}_2=\bmat{0& 0\\-\frac{1}{\sqrt{2}}& 1}}.
\end{equation*}
\end{Example}
\noindent Notice that
\begin{equation*}
\mathrm{S}_1^{\ell}=\bmat{1&\frac{\ell}{\sqrt{2}}\\0 & 1},\
\qquad \mathrm{S}_2^m=\mathrm{S}_2.
\end{equation*}
and
\begin{equation*}
\mathrm{S}_1^{\ell}\mathrm{S}_2^m=\bmat{\frac{\ell}{\sqrt{2}}\\1}\bmat{-\frac{1}{\sqrt{2}}& 1}.
\end{equation*}
Thus we have
\begin{equation*}
\rho(\mathrm{S}_1^{\ell}\mathrm{S}_2^m)=|1-\frac{\ell}{2}|.
\end{equation*}
Hence
\begin{equation*}\begin{split}
\rho(\bS)&=\max_{\ell, m\in \N}\sqrt[\ell+m]{|1-\frac{\ell}{2}|}=\sqrt[11]{4}\\
&\approx1.134313> \max\left\{\rho(\mathrm{S}_1),\rho(\mathrm{S}_2)\right\},
\end{split}\end{equation*}
where the maximum is attained at $(\ell,m)=(10,1)$ with the optimal sequence $\mathrm{S}_1^{10}\mathrm{S}_2$.

\begin{Example}\label{exmple4}
\begin{equation*}
\bS=\set{\mathrm{S}_1=\bmat{1& \eps&0&0\\0& 1&\eps&0\\0&0&1&\eps\\0&0&0&1},
\mathrm{S}_2=\bmat{1& -1&0&1\\1& -1&0&1\\1& -1&0&1\\1& -1&0&1}},
\end{equation*}
where $\eps>0$ is a parameter.
\end{Example}

\noindent Notice that
\begin{equation*}
\mathrm{S}_1^{\ell}=\left[
\begin{array}{llll}
 1 & \ell  \eps  & \frac{1}{2} (\ell-1 ) \ell  \eps^2 & \frac{1}{6} (\ell-2) (\ell-1 ) \ell  \eps^3 \\
 0 & 1 & \ell  \eps  & \frac{1}{2} (\ell-1) \ell  \eps^2 \\
 0 & 0 & 1 & \ell  \eps  \\
 0 & 0 & 0 & 1
\end{array}
\right],\
\qquad \mathrm{S}_2^m=\mathrm{S}_2.
\end{equation*}
and
\begin{equation*}
\mathrm{S}_1^{\ell}\mathrm{S}_2^m=\bmat{
1 + \ell  \eps + \frac{1}{2} (\ell-1 ) \ell  \eps^2 + \frac{1}{6} (\ell-2) (\ell-1 ) \ell  \eps^3 \\
 1 + \ell  \eps + \frac{1}{2} (\ell-1) \ell  \eps^2 \\
 1 + \ell  \eps  \\
 1
}
\bmat{1&-1&0&1}.
\end{equation*}
Thus we have
\[
\rho(\mathrm{S}_1^{\ell}\mathrm{S}_2^m)=\frac{1}{6} (\ell-2) (\ell-1 ) \ell  \eps^3+1.
\]
Hence
\begin{equation*}\begin{split}
\rho(\bS)&=\max_{\ell, m\in \N}\sqrt[\ell+m]{
\frac{1}{6} (\ell-2) (\ell-1 ) \ell  \eps^3+1}\\
&=\max_{\ell\ge 3}\sqrt[\ell+1]{\frac{1}{6} (\ell-2) (\ell-1 ) \ell  \eps^3+1}\\
&=\sqrt[\ell_{\eps}+1]{\frac{1}{6} (\ell_{\eps}-2) (\ell_{\eps}-1 ) \ell_{\eps}  \eps^3+1},
\end{split}\end{equation*}
where the maximum is assumed to be achieved at $\ell=\ell_{\eps}$.
One can show that $\ell_{\eps}\to \infty$ as $\eps\to 0$. Numerical experiments indicate that  value of $\ell_{\eps}$  increases very quickly with respect to $\frac{1}{\eps}$.
Thus for any given integer $L>0$, one always can find a corresponding constant $\eps>0$, such that $\sqrt[L]{\rho(\mathrm{S}_{i_1}\mathrm{S}_{i_2}\cdots \mathrm{S}_{i_L})}<\rho(\bS)$.

This argument also can be easily shown by the following two-dimensional example.

\begin{Example}\label{example5}
\begin{equation*}
\bS=\set{ S_1=\bmat{1& \epsilon \\0& 1},\mathrm{S}_2=\bmat{1& -1\\1& -1}}.
\end{equation*}
\end{Example}

\noindent Since
\begin{equation*}
\mathrm{S}_1^{\ell}=\bmat{1&\ell \epsilon\\ 0& 1}, \ \qquad \mathrm{S}_2^2=0
\end{equation*}
we have
\begin{equation*}
\mathrm{S}_1^{\ell}\mathrm{S}_2 =\bmat{1+\ell \epsilon\\1}\bmat{1& -1}.
\end{equation*}
Then
\begin{equation*}
\rho(\mathrm{S}_1^{\ell}\mathrm{S}_2)=\ell \epsilon.
\end{equation*}
Hence
\begin{equation*}
\rho(\bS)
=\max_{\ell\in \N} \sqrt[\ell+1]{\ell \epsilon}.
\end{equation*}
Given any specified length $L$, let $\epsilon=\frac{1}{L+1}$, then
\begin{equation*}\begin{split}
\rho(\bS)&=\max_{\ell\in \N} \sqrt[\ell+1]{\frac{\ell}{L+1}}\\
&\ge1\\
&>\max_{1\le\ell\le L} \sqrt[\ell+1]{\frac{\ell}{L+1}},
\end{split}\end{equation*}
where the last strictly inequality implies that, for the chosen $\epsilon=\frac{1}{L+1}$,
the intended optimal sequence will never be found within the  length $L$.
This special example represents the challenge even if we know the spectral finiteness property holds.
Therefore, any algorithms depending on the search of the length of optimal sequence will suffer from a high computational cost.

\section{Concluding remarks}\label{sec4}
In this paper, we have proved that the spectral finiteness property holds for every pairs of real $d\times d$ matrices $\mathrm{S}_1, \mathrm{S}_2$, if one of $\mathrm{S}_1,\mathrm{S}_2$ has rank $1$; see Theorem~\ref{thm1.3}. Under our context, $\mathrm{S}_1$ and $\mathrm{S}_2$ might be neither symmetric, nor commutative, and nor rational. In addition, our argument does not involve any polytope norms.

Recall that a matrix $A=[a_{ij}]$ is called a binary matrix (resp. sign-matrix), provided that every entries $a_{ij}$ belong to $\{0,1\}$ (resp. $\{-1,0,1\}$).
In \cite[Theorem~4]{JB08}, R.~Jungers and V.~Blondel proved that
\begin{itemize}
\item The finiteness property holds for all sets of nonnegative rational square matrices if and only it holds for all pairs of binary square matrices.

\item The finiteness property holds for all sets of rational square matrices if and only it holds for all pairs of square sign-matrices.
\end{itemize}
Moreover, the following two positive results are already known.
\begin{itemize}
\item The finiteness property holds for every pairs of $2\times 2$ binary matrices (\cite{JB08}).

\item The finiteness property holds for every pairs of $2\times 2$ sign-matrices (\cite{CGSZ10}).
\end{itemize}
However, the above two results cannot imply the finiteness property for every pairs of $2\times 2$ (nonnegative) rational matrices; this is because when one reduces, following the framework of \cite[Theorems~2 and 4]{JB08}, a pair of $2\times 2$ rational matrices $\mathrm{S}_1, \mathrm{S}_2$ to a pair of binary or sign-matrices $\widetilde{\mathrm{S}}_1,\widetilde{\mathrm{S}}_2$, the size of $\widetilde{\mathrm{S}}_1,\widetilde{\mathrm{S}}_2$ becomes $2m\times 2m$ and in general $m$ would be sufficiently large.

So, our result is essentially new and our approach has the additional interest.


\begin{thebibliography}{99}
\bibitem{Bar}
   \newblock {\sc N. Barabanov},
   \newblock {\it Lyapunov indicators of discrete inclusions I--III},
   \newblock {Autom. Remote Control, 49 (1988), pp.~152--157, 283--287, 558--565}.

\bibitem{Bar89}
   \newblock {\sc N.\,E.~Barabanov},
   \newblock {\it An absolute characteristic exponent of a class of linear nonstationary systems of differential equations},
   \newblock {Siberian Math. J. 29 (1989)~521--530}.

\bibitem{Bar05}
   \newblock {\sc N. Barabanov},
   \newblock {\it Lyapunov exponent and joint spectral radius: some known and new results},
   \newblock {in: Proc. 44th IEEE Conf. Decision and Control and Eur. Control Conf., Seville, Spain, 2005, pp.~2332--2336}.

\bibitem{BW92}
  \newblock {\sc M.A.~Berger and Y.~Wang},
  \newblock {\it Bounded semigroups of matrices},
  \newblock {Linear Algebra Appl., {166} (1992), pp.~21--27}.

\bibitem{BJP06}
  \newblock {\sc V.D.~Blondel, R.~Jungers, and V.Yu.~Protasov},
  \newblock {\it On the complexity of computing the capacity of codes that avoid forbidden difference patterns},
  \newblock {IEEE Trans. Inform. Theory, 52 (2006), pp.~5122--5127}.

\bibitem{BN05}
  \newblock {\sc V.D.~Blondel and Y.~Nesterov},
  \newblock {\it Computationally efficient approximations of the joint spectral radius},
  \newblock {SIAM J. Matrix Anal. Appl., {27} (2005), pp.~256--272}.

\bibitem{BTV}
  \newblock {\sc V.D.~Blondel, J.~Theys, and A.A.~Vladimirov},
  \newblock {\it An elementary counterexample to the finiteness conjecture},
  \newblock {SIAM J. Matrix Anal. Appl., {24} (2003), pp.~963--970}.

\bibitem{Bochi}
  \newblock {\sc J.~Bochi},
  \newblock {\it Inequalities for numerical invariants of sets of matrices},
  \newblock {Linear Algebra Appl., {368} (2003), pp.~71--81}.

\bibitem{BM}
  \newblock {\sc T.~Bousch, and J.~Mairesse},
  \newblock {\it Asymptotic height optimization for topical IFS, Tetris heaps and the finiteness conjecture},
  \newblock {J. Amer. Math. Soc., {15} (2002), pp.~77--111}.

\bibitem{CZ00}
   \newblock {\sc Q.~Chen} and {\sc X.~Zhou},
   \newblock {\it Characterization of joint spectral radius via trace},
   \newblock {Linear Algebra Appl., {315} (2000), pp.~175--188}.

\bibitem{CGSZ10}
  \newblock {\sc A.~Cicone, N.~Guglielmi, S.~Serra-Capizzano} and {\sc M.~Zennaro},
  \newblock {\it Finiteness property of pairs of $2\times 2$ sign-matrices via real extremal polytope norms},
  \newblock {Linear Algebra Appl., 432 (2010), pp.~796--816}.

\bibitem{Dai-JMAA}
    \newblock {\sc X.~Dai},
    \newblock {\it Extremal and Barabanov semi-norms of a semigroup generated by a bounded family of matrices},
    \newblock {J. Math. Anal. Appl., 379 (2011), pp.~827--833}.

\bibitem{Dai-JDE}
       \newblock {\sc X.~Dai},
       \newblock {\it Weakly Birkhoff recurrent switching signals, almost sure and partial stability of linear switched dynamical systems},
       \newblock {J. Differential Equations, 250 (2011), pp.~3584--3629}.

\bibitem{Dai-sicon}
       \newblock {\sc X.~Dai},
       \newblock {\it Criterion of stabilizability for switching systems with solvable linear approximations}, Preprint 2010.

\bibitem{DHX11}
       \newblock {\sc X.~Dai, Y.~Huang, and M.~Xiao},
       \newblock {\it Periodically switched stability induces exponential stability of discrete-time linear switched systems in the sense of Markovian probabilities},
       \newblock {Automatica J. IFAC, 47 (2011), pp.~1512--1519}.

\bibitem{DHX-ERA}
  \newblock {\sc X.~Dai, Y.~Huang, and M.~Xiao},
  \newblock {\it Realization of joint spectral radius via ergodic theory},
  \newblock {Electron. Res. Announc. Math. Sci., {18} (2011), 22--30}.


\bibitem{DL92-01}
  \newblock {\sc I.~Daubechies and J.C.~Lagarias},
  \newblock {\it Sets of matrices all infinite products of which converge},
  \newblock {Linear Algebra Appl., {161} (1992), pp.~227--263. Corrigendum/addendum, {327} (2001), pp.~69--83}.

\bibitem{DL92}
  \newblock {\sc I.~Daubechies and J.C.~Lagarias},
  \newblock {\it Two-scale difference equations. II. Local regularity, infinite products of matrices and fractals},
  \newblock {SIAM J. Math. Anal., {23} (1992), pp.~1031--1079}.

\bibitem{DST}
  \newblock {\sc J.M.~Dumont, N.~Sidorov, and A.~Thomas},
  \newblock {\it Number of representations related to a linear recurrent basis},
  \newblock {Acta Arith., 88 (1999), pp.~371--396}.

\bibitem{El}
  \newblock {\sc L.~Elsner},
  \newblock {\it The generalized spectral-radius theorem: An analytic-geometric proof},
  \newblock {Linear Algebra Appl., {220} (1995), pp.~151--159}.

\bibitem{Fen}
  \newblock {\sc N.~Fenichel},
  \newblock {\it Persistence and smoothness of invariant manifolds for flows},
  \newblock {Indiana Univ. Math. J., 21 (1971), pp.~193--226}.

\bibitem{Gri96}
  \newblock {\sc G.~Gripenberg},
  \newblock {\it Computing the joint spectral radius},
  \newblock {Linear Algebra Appl., {234} (1996), pp.~43--60}.

\bibitem{GWZ05}
  \newblock {\sc N.~Guglielmi, F.~Wirth, and M.~Zennaro},
  \newblock {\it Complex polytope extremality results for families of matrices},
  \newblock {SIAM J. Matrix Anal. Appl., 27 (2005), pp.~721--743}.

\bibitem{GZ}
  \newblock {\sc N.~Guglielmi and M.~Zennaro},
  \newblock {\it On the zero-stability of a family of variable stepsize multistep methods: The spectral radius approach},
  \newblock {Numer. Math., 88 (2001), pp.~445--458}.

\bibitem{GZ01}
  \newblock {\sc N.~Guglielmi and M.~Zennaro},
  \newblock {\it On the asymptotic properties of a family of matrices},
  \newblock {Linear Algebra Appl., 322 (2001), pp.~169--192}.

\bibitem{GZ08}
  \newblock {\sc N.~Guglielmi and M.~Zennaro},
  \newblock {\it An algorithm for finding extremal complex polytope norms for matrix families},
  \newblock {Linear Algebra Appl., 428 (2008), pp.~2265--2282}.

\bibitem{GZ09}
  \newblock {\sc N.~Guglielmi and M.~Zennaro},
  \newblock {\it Finding extremal complex polytope norms for a family of real matrices},
  \newblock {SIAM J. Matrix Anal. Appl., 31 (2009), pp.~602--620}.

\bibitem{Gur95}
  \newblock {\sc L.~Gurvits},
  \newblock {\it Stability of discrete linear inclusions},
  \newblock {Linear Algebra Appl., {231} (1995), pp. 47--85}.

\bibitem{HMST}
  \newblock {\sc K.G.~Hare, I.D.~Morris, N.~Sidorov, and J.~Theys},
  \newblock {\it An explicit counterexample to the Lagarias-Wang finiteness conjecture},
  \newblock {Adv. Math., 226 (2011), pp.~4667--4701}.

\bibitem{HS95}
  \newblock {\sc C.~Heil and G.~Strang},
  \newblock {\it Continuity of the joint spectral radius: application to wavelets},
  \newblock {Linear Algebra for Signal Processing, The IMA Volumes in Mathematics and its Applications, vol. 69, Springer, New York, 1995, pp.~51--61}.

\bibitem{Jun09}
  \newblock {\sc R.~Jungers},
  \newblock {\it The joint spectral radius, theory and applications},
  \newblock {Lecture Notes in Control and Information Sciences, vol. 385, Springer-Verlag, Berlin, 2009}.

\bibitem{JB08}
  \newblock {\sc R.~Jungers and V.~Blondel},
  \newblock {\it On the finiteness properties for rational matrices},
  \newblock {Linear Algebra Appl., 428 (2008), pp.~2283--2295}.

\bibitem{JP09}
  \newblock {\sc R.~Jungers and V.Y.~Protasov},
  \newblock {\it Counterexamples to the complex polytope extremality conjecture},
  \newblock {SIAM J. Matrix Anal. Appl. 31 (2009), pp.~404--409}.

\bibitem{JPB08}
  \newblock {\sc R.~Jungers, V.Y.~Protasov, and V.~Blondel},
  \newblock {\it Efficient algorithms for deciding the type of growth of products of integer matrices},
  \newblock {Linear Algebra Appl. 428 (2008), pp.~2296--2311}.

\bibitem{Koz90}
  \newblock {\sc V.S.~Kozyakin},
  \newblock {\it Algebraic unsolvability of a problem on the absolute stability of desynchronized systems},
  \newblock {Autom. Remote Control, 51 (1990), pp.~754--759}.

\bibitem{Koz05}
  \newblock {\sc V.S.~Kozyakin},
  \newblock {\it Extremal norms, discontinuous circle maps and a counterexample to the finiteness conjecture},
  \newblock {Information Processes, 5 (2005), pp.~301--335}.

\bibitem{Koz07}
   \newblock {\sc V.S.~Kozyakin},
   \newblock {\it Structure of extremal trajectories of discrete linear systems and the finiteness conjecture},
   \newblock {Autom. Remote Control, {68} (2007), pp.~174--209}.

\bibitem{Koz09-1}
  \newblock {\sc V.S.~Kozyakin},
  \newblock {On the computational aspects of the theory of joint spectral radius},
  \newblock {Dokl. Math., 80 (2009), pp.~487--491}.

\bibitem{Koz09-2}
  \newblock {\sc V.S.~Kozyakin},
  \newblock {\it On accuracy of approximation of the spectral radius by the gelfand formula},
  \newblock {Linear Algebra Appl., 431 (2009), PP.~2134--2141}.

\bibitem{Koz10-DEDS}
  \newblock {\sc V.S.~Kozyakin},
  \newblock {\it On explicit a priori estimates of the joint spectral radius by the generalized gelfand formula},
  \newblock {Differ. Equ. Dyn. Syst., 18 (2010), pp.~91--103}.

\bibitem{Koz10}
   \newblock {\sc V.S.~Kozyakin},
   \newblock {\it Iterative building of Barabanov norms and computation of the joint spectral radius for matrix sets},
   \newblock {Discrete Contin. Dyn. Syst. Ser. B, 14 (2010), pp.~143--158}.

\bibitem{LW95}
  \newblock {\sc J.C.~Lagarias and Y.~Wang},
  \newblock {\it The finiteness conjecture for the generalized spectral radius of a set of matrices},
  \newblock {Linear Algebra Appl., 214 (1995), pp.~17--42}.

\bibitem{LD06-1}
  \newblock {\sc J.W.~Lee and G.E.~Dullerud},
  \newblock {\it Uniform stabilization of discrete-time switched and Markovian jump linear systems},
  \newblock {Automatica J. IFAC, 42 (2006), pp.~205--218}.

\bibitem{LD06-2}
  \newblock {\sc J.W.~Lee and G.E.~Dullerud},
  \newblock {\it Optimal disturbance attenuation for discrete-time switched and Markovian jump linear systems},
  \newblock {SIAM J. Control Optim., 45 (2006), pp.~1329--1358}.

\bibitem{LHM}
  \newblock {\sc D.~Liberzon, J.P.~Hespanha, and A.S.~Morse},
  \newblock {\it Stability of switched systems: a Lie-algebraic condition},
  \newblock {Syst. Control Lett., 37 (1999), pp.~117--122}.

\bibitem{Mae96}
  \newblock {\sc M.~Maesumi},
  \newblock {\it An efficient lower bound for the generalized spectral radius of a set of matrices},
  \newblock {Linear Algebra Appl., 240 (1996), pp.~1--7}.

\bibitem{Mae98}
  \newblock {\sc M.~Maesumi},
  \newblock {\it Calculating joint spectral radius of matrices and H\"{o}lder exponent of wavelets},
  \newblock {Approximation theory IX, Vol. {2} (Nashville, TN, 1998), Innov. Appl. Math., Vanderbilt Univ. Press, Nashville, TN, 1998, pp.~205--212}.

\bibitem{Morris10}
  \newblock {\sc I.D.~Morris},
  \newblock {\it Criterion for the stability of the finiteness property and for the uniquenss of Barabanov norms},
  \newblock {Linear Algebra Appl., 433 (2010), pp.~1301--1311}.

\bibitem{Morris-am}
  \newblock {\sc I.D.~Morris},
  \newblock {\it A rapidly-converging lower bound for the joint spectral radius via multiplicative ergodic theory},
  \newblock {Adv. Math., 225 (2010), pp.~3425--3445}.

\bibitem{MOS}
  \newblock {\sc B.E.~Moision, A.~Orlitsky, and P.H.~Siegel},
  \newblock {\it On codes that avoid specified differences},
  \newblock {IEEE Trans. Inform. Theory, {47} (2001), pp.~433--442}.

\bibitem{PJ08}
  \newblock {\sc P.A.~Parrilo and A.~Jadbabaie},
  \newblock {\it Approximation of the joint spectral radius using sum of squares},
  \newblock {Linear Algebra Appl., 428 (2008), pp.~2385--2402}.

\bibitem{Pro96}
  \newblock {\sc V.Y.~Protasov},
  \newblock {\it The joint spectral radius and invariant sets of linear operators},
  \newblock {Fundam. Prikl. Mat., {2} (1996), pp.~205--231}.

\bibitem{Pro97}
  \newblock {\sc V.Y.~Protasov},
  \newblock {\it A generalized joint spectral radius: A geometric approach},
  \newblock {Izv. Math., {61} (1997), pp.~995--1030}.

\bibitem{PR91}
  \newblock {\sc E.S.~Pyatnitski\u{i} and L.B.~Rapoport},
  \newblock {\it Periodic motion and tests for absolute stability on nonlinear nonstationary systems},
  \newblock {Autom. Remote Control, {52} (1991), pp.~1379--1387}.

\bibitem{RS}
  \newblock {\sc G.-C.~Rota and G.~Strang},
  \newblock {\it A note on the joint spectral radius},
  \newblock {Indag. Math., {22} (1960), pp.~379--381}.

\bibitem{SWP}
   \newblock {\sc M.-H.~Shih, J.-W.~Wu, and C.-T.~Pang},
   \newblock {\it Asymptotic stability and generalized Gelfand spectral radius formula},
   \newblock {Linear Algebra Appl., {252} (1997), pp.~61--70}.

\bibitem{SWMW07}
   \newblock {\sc R.~Shorten, F.~Wirth, O.~Mason, K.~Wulff, and C.~King},
   \newblock {\it Stability criteria for switched and hybrid systems},
   \newblock {SIAM Rev., {49} (2007), pp. 545--592}.

\bibitem{Theys}
  \newblock {\sc J.~Theys},
  \newblock {\it Joint Spectral Radius: Theory and Approximations},
  \newblock {Ph.D. thesis, Universit\'{e} Catholique de Louvain, 2005}.

\bibitem{TB97}
  \newblock {\sc J.N.~Tsitsiklis and V.D.~Blondel},
  \newblock {\it The Lyapunov exponent and joint spectral radius of pairs of matrices are hard\,--\,when not impossible\,--\,to compute and to approximate},
  \newblock {Math. Control Signals Systems, {10} (1997), pp.~31--40}.

\bibitem{Wirth02}
  \newblock {\sc F.~Wirth},
  \newblock {\it The generalized spectral radius and extremal norms},
  \newblock {Linear Algebra Appl., 342 (2002), pp.~17--40}.

\bibitem{Wirth05}
  \newblock {\sc F.~Wirth},
  \newblock {\it The generalized spectral radius is strictly increasing},
  \newblock {Linear Algebra Appl., 395 (2005), pp.~141--153}.
\end{thebibliography}
\end{document}